\documentclass{article} \usepackage{amssymb,amsmath,amsthm,color}
\usepackage[latin1]{inputenc}
\usepackage{hyperref}\usepackage[numbers]{natbib}\usepackage{url}
\newif
\ifpdf\ifx\pdfoutput\undefinedles
\pdffalse\else
\pdfoutput=1\pdftrue\fi\ifpdf\pdfcompresslevel=9\fi
\makeatletter
\renewcommand\section{
\@startsection {section}{1}{\z@}%
{-3.5ex \@plus -1ex \@minus -.2ex}%
{2.3ex \@plus.2ex}%
{\bf\centering\normalsize}}
\renewcommand\subsection{\@startsection{subsection}{2}{\z@}%
  {-3.25ex\@plus -1ex \@minus -.2ex}%
{1.5ex \@plus .2ex}%
{\normalfont\normalsize\bfseries}}
\renewcommand\subsubsection{\@startsection{subsubsection}{3}{\z@}%
{-3.25ex\@plus -1ex \@minus -.2ex}%
{1.5ex \@plus .2ex}%
{\normalfont\normalsize\bfseries}}
\renewcommand\paragraph{\@startsection{paragraph}{4}{\z@}%
{3.25ex \@plus1ex \@minus.2ex}%
{-1em}%
{\normalfont\normalsize\bfseries}}
\renewcommand\subparagraph{\@startsection{subparagraph}{5}{\parindent}%
{3.25ex \@plus1ex \@minus .2ex}%
{-1em}%
{\normalfont\normalsize\bfseries}} \makeatother
\newcommand{\TheoremText}{Theorem}
\newcommand{\ConditionText}{Condition}
\newcommand{\PropositionText}{Proposition}
\newcommand{\LemmaText}{Lemma} \newcommand{\CorollaryText}{Corollary}
\newcommand{\ProofText}{Proof} \newcommand{\AxiomText}{Axiom}
\newcommand{\DefinitionText}{Definition}
\newcommand{\RemarkText}{Remark} \newcommand{\NoteText}{Note}
\newcommand{\ExampleText}{Example}
\newcommand{\ConventionText}{Convention}
\newcommand{\ExerciseText}{Exercise}
\newcommand{\WarningText}{Warning} \newcommand{\ProblemText}{Problem}

\newtheorem{theorem}{\TheoremText}
\newtheorem{condition}{\ConditionText}
\newtheorem{proposition}{\PropositionText}

\newtheorem{definition}{\DefinitionText}
\newtheorem{problem}{\ProblemText} \newtheorem{varremark}{\RemarkText}
\newtheorem{varnote}{\NoteText} \newtheorem{varexample}{\ExampleText}
\newtheorem{varconvention}{\ConventionText}
\newtheorem{varwarning}{\WarningText}
\newenvironment{remark}{\begin{varremark}\em}{\em\end{varremark}}

\newenvironment{example}{\begin{varexample}\em}{\em\end{varexample}}

\renewenvironment{proof}{
  \noindent\textbf{\ProofText}\ }{\hspace*{\fill}
  \begin{math}\Box\end{math}\medskip}
\newenvironment{proof*}[1]{
  \noindent\textbf{#1\ }}{\hspace*{\fill}
  \begin{math}\Box\end{math}\medskip}
\newcounter{exercisenr} 

\title{State Dependent Utility\thanks{I want to thank professor M. M.
    Rao for a reading of the first version of this paper and
    suggestions made on it that led to an important improvement.  All
    remaining errors are mine. This work was 
    supported by COLCIENCIAS (Colombian national science foundation)
    with grant number 284-2003, Universidad EAFIT, and Universidad
    Nacional de Bogotá.  }}
\author{Jaime A. Londoño\\
  Departamento de Matemáticas\\
  Universidad Nacional de Colombia, Bogotá, Colombia\\
  jlondono@math.ucr.edu}
\begin{document}
\maketitle
\tolerance=200 \setlength{\emergencystretch}{2em}

\begin{abstract}
  We propose a new approach to utilities that is consistent with
  state-dependent utilities.  In our model utilities reflect the
  level of  consumption satisfaction of flows of cash in future
  times as they are valued when the economic agents are
  making their consumption and investment decisions.  The theoretical
  framework used for the model is one proposed by the author in
  \emph{Dynamic State Tameness (arXiv:math.PR/0509139)}.  The proposed framework  is
a generalization of the theory of Brownian
flows and can be applied to those processes that are the
solutions of classical Itô stochastic differential equations, even
when the volatilities and drifts are just locally $\delta$-Hölder
continuous for some $\delta>0$.    We develop the martingale
methodology for the solution of the problem of optimal consumption and
investment.  Complete solutions of the optimal consumption and portfolio problem are
obtained in a very general setting which includes several  functional
forms for utilities in the  current literature, and consider general
restrictions on minimal wealths.  As a secondary result we obtain a
suitable representation for straightforward numerical computations of
the optimal consumption and investment strategies. 
\end{abstract}
\section[Introduction]{Introduction}
\label{sec:introduction}
The problem of optimal consumption and investment for a ``small investor'' whose
actions do not influence market prices is at the core of portfolio
management and it is the building block for the development of equilibrium theory.  The modern treatment of this problem when asset prices
follows Itô processes started with the seminal works 
of \citet{Merton69} and \citet{Merton71}. Using  a ``martingale''
approach,  \citet{Cox_etal89}, and \citet{Karatzas_etal87} solved the
problem in more general settings in the case of complete markets. A
representation formula is derived in \citet{Ocone_Karatzas91} in terms
of expectation of random variables which involve Malliavin derivatives
of the coefficients of the model.  The latter gives theoretical formulas for
optimal portfolios and consumption strategies.

In order to obtain numerical representations for the structure of the
optimal portfolios and consumption processes it is  natural to
use methods based on the dynamic-programming approach.
However numerical schemes based on PDEs become increasingly difficult to
evaluate when the dimension of the underlying state variable increase
and even standard techniques are somehow inappropriate for the solution
of the PDE that arises in small dimensions (\citet{Dangl_Wirl04}).  As a
result attention has been directed to models admitting closed form
solutions  (i.e \citet{Wachter02}, \citet{Kim_Omberg96}, \citet{Abramham_Poncet01}), specifications which are computationally tractable based on
dynamic programming techniques
(\citet{Brennan_etal97}, \citet{Brennan98}, \citet{Brennan_Xia01}, \citet{Campbell_etal04} ), discrete time models based on
approximated Euler equations (\citet{Balduzzi_Lynch99},
\citet{Dammon_etal01}, \citet{Campbell_Viceira99}, and
\citet{Campbell_Viceira99b}) or Monte Carlo tecniques
(\citet{Cvitanic_etal03} and \citet{Detemple_etal03}).

However, the main drawback of the standard  models for optimal consumption and
investment is their lack of agreement with empirical data. These
inconsistencies are documented with the name of several puzzles such
as the ``equity premium puzzle'' (\citet{Mehra_and_Prescott85}),
the ``risk-free rate puzzle'' (\citet{Weil89}), and ``risk-aversion
puzzle'' (\citet{Jackwerth00}) .  In order to address
these problems several generalizations  have been proposed.  One of these
models habit formation for consumers. Some examples are
\citet{Constantinides90}, \citet{Hindy_and_Huang93},
\citet{Hindy_etal97}, \citet{Hindy_etal97b}.  Another approach is the
construction of  recursive utility.   Some references for this are  \citet{Duffie_and_Epstein92a},
\citet{Epstein_and_Zin89} and \citet{Lazrak03}.   A different way to
account for the discrepancy of theory and  empirical data is the
assumption of transaction costs for changes in consumption levels.
Some references for this approach are \citet{Magill_Constantinides76},
\citet{Shreve_Soner94}, \citet{Davis_Norman90} to cite a few.

In fact one of the reasons why the standard utility models fail
to fit economic behavior might be the fact that state independent
utilities are not appropriate for modeling the behavior of human
beings.  For instance, see \citet{Karni93}, \citet{Karni93b}.  Partly
motivated by the above, some literature in finance has focused on
state-dependent utilities to explain the behavior of individual
consumers and investors and of financial variables. Some recent
references are \citet{Chabi-yo_etal05}, \citet{Melino03} and
\citet{Danthine_etal04} among others.

In this paper we propose a new approach for utilities.  Mathematically
it looks similar to the standard model for utilities, but we
interpreted it in a way that is consistent with state-dependent
utilities.  The traditional approach is to consider that utilities  reflect the
level of ``happiness'' for   consumption levels in the future (discounted by the value of money in a bank account). See \citet{Karatzas98}.  We
believe that in complete markets where an agent can hedge any flow of
money, this view is inappropriate.  To simplify things let us assume that
a time horizon $[0,T]$ is fixed, and $0$ is the time when  the agent is making his 
decision on an optimal consumption level  and investment and there is
not any preference for terminal  wealth.  \emph{Our guiding principle
  is the  believe that agents have
utilities for consumption of flows of money in future times  as they
are valued (by the market)  at the time when they are making their consumption and
investment decisions.}  Another way to look at this,  is that
people tend to value things according to their social and economic
context, instead of just looking at quantitative values.  For instance
people tend to appreciate more  the ability to have  enough
money to pay off their debts in depression times than the ability to
buy  luxuries in good times. 
The above remark  changes completely the optimization problem to
consider, with the advantage that most of the tools used to solve the
old problem can be used in this setting. In particular the martingale
methodology is available.  A consequence of adopting this approach
is that complete solutions of the optimal consumption and portfolio problem are
obtained in a very general setting that includes several of the functional
forms for utilities in the literature, and considers quite general
restrictions of minimal wealth.  As a secondary result we obtain
suitable representation for straightforward numerical computations of
the optimal consumption and \emph{investment} strategies.

The theoretical framework used to solve the above problem is the one
proposed in \citet{Londono2005a}.  In this introduction we just named the
processes described in the cited paper as consistent measurable
processes; these are processes whose evolution between any two times only
depends on \emph{the evolution} of the underlying Brownian motion and
satisfies some consistency conditions.  The methodology described can be used in processes that are
a generalization of Brownian
flows (\citet{Kunita90}), and is applicable to those processes that are the
solutions of classical Itô stochastic differential equations, even
when the volatilities and drifts are just locally $\delta$-Hölder
continuous for some $\delta>0$.  It  can also  be straightforwardly adapted to
stochastic volatility models whose underlying drifts
and volatilities evolution are described by classical stochastic Itô
differential equations. 
As discussed in \citet{Londono2005a} the
theoretical framework is potentially useful when the underlying
randomness is generated by a (not necessarily continuous) Lévy
process. At the same time the theoretical framework described in
\citet{Londono2005a} is a particular case  of the theory of arbitrage and
valuation  presented in \citet{Londono04}.  To the best of our
knowledge this theory of arbitrage and valuation is the most general
existing theory in the case of (continuous) semimartingales driven by
Brownian filtrations with continuous coefficients. 

A brief explanation  of the heuristics of using consistent measurable
processes in finance  (in the problem of optimal consumption and
investment) is given below.   
When an agent enters an economy, it would be desirable that any
decision made during his lifetime would depend only on events that
are happening within the time framework defined by the current time
and the entering time, and on the state of the economy when he enters;
we hope to be able to show to the reader that this is the case
throughout this paper.  Although this might be a simplification that
is likely not to be truthful, in the long run events in the past
(before the agent entered the economy) become neglectful.
Even more important, past aware models necessarily make the entering
time a  privileged
one; the entering time would be the only time when no past
information is required (besides the information contained in the
current state of the model). Moreover, it would imply that the
decision making process of agents would depend on their entering time,
regardless of identical preferences over their lifetime. The previous
remark also applies to the standard basic time-additive utility
function maximization model ~(\citet{Karatzas98}). Let us expand on
this point further.  Assume that two agents decide to start investing
in the market at two different times, let's say the first agent starts
investing at time $t=0$, and the second agent starts investing at time
$t=1$; assume that both agents have the same live expectation (to
simplify things, assume that both are planing to leave the economy at
time $t=2$), and they do not have any additional source of income.
Moreover, assume that both agents have identical preferences over
consumption partners for times $1\leq t\leq 2$, and at time $t=1$ they
have identical amounts of wealth.  We believe that any reasonable model
should produce identical levels of consumption and investment,
assuming that they are the result of some kind of optimization
procedure.  Unfortunately, this is not the case for the basic
time-additive utility function maximization model
(\citet{Karatzas98}) when the processes are allowded to be general
semimartingales.
 In fact, when an agent chooses his optimal
consumption process at time $0$ over the time framework $[0,2]$, the
agent is free to use portfolios and consumption processes that at any
time ``know'' all the previous information between the given time and
the initial time ($t=0$). In particular the consumption and portfolio
processes for times $t\in[1,2]$ would be ``aware'' of all the
information before time $t=1$, and therefore \emph{a priori} they
would be different from the ones that are the result of the
optimization procedure that the second agent is making at time $t=1$
(with no previous information gathered).  Note that the above does not
hold for models which are solutions of stochastic differential
equations whose coefficients are deterministic functions.

Next we describe the contents of the paper.  In section
\ref{sec:model} we review the model and definitions presented in
\citet{Londono2005a}, and review the definitions of utility that we use
in this paper.  In section \ref{sec:martingale-approach} we present a
martingale methodology needed to address the cited problem for  the model described in
\citet{Londono2005a}, that plays the role of the martingale methodology of
\citet{Cox_etal89}, \citet{Karatzas_etal87} and
\citet{Ocone_Karatzas91} in the current context.  Finally, in section
\ref{sec:cons-portf-optim} we present the main results on optimal
consumption and investment.

\section[The model]{The model}
\label{sec:model}
First we introduce some notation which will be frequently used in this
paper.  Let $\mathbb{D}\subset\mathbb{R}^k$ be a open connected set.
Let $m$ be a non-negative integer.  We denote by
$C^{m,\delta}(\mathbb{D}\colon\mathbb{R}^n)$ the the Fréchet space of
$m$-times continuous differentiable functions whose $m$-order
derivatives are $\delta$-H\"older continuous with semi-norms
$\|f\|_{m,\delta\colon K}$ defined in \citet[Section 3.1]{hK90} where
$K\subset \mathbb{D}$ is a compact set and $0\leq\delta\leq 1$.  In
case $m=0$ (or $\delta=0$) we  denote
$C^{m,\delta}(\mathbb{D}\colon\mathbb{R}^n)$ simply by
$C^{\delta}(\mathbb{D}\colon\mathbb{R}^n)$
($C^m(\mathbb{D}\colon\mathbb{R}^n)$). 

We assume a $d$-dimensional Brownian Motion $\{W(t),\mathcal{F}_t;
0\leq t\leq T\}$ starting at $0$ and  defined on a complete probability
space $(\Omega,\mathcal{F},\mathbf{P})$ where
$\mathcal{F}=\mathcal{F}_T$ and $\{\mathcal{F}_t, 0\leq t\leq T\}$ is
the $\mathbf{P}$ augmentation by the null sets of the natural
filtration $\mathcal{F}^{W}_t=\sigma(W(s),0\leq s\leq t)$.  Let
$(\mathcal{F}_{s,t})=\left\{\mathcal{F}_{s,t}, 0\leq s\leq t
  \leq T\right\} $ be the two parameter filtration where
$\mathcal{F}_{s,t}$ is the smallest sub $\sigma$-field containing all
null sets and $\sigma(W_s(u)\mid s\leq u\leq t)$, where $W_s(u)\equiv
W(u)-W(s)$.   For each $0\leq
s\leq T$ we also define the  $\sigma$-field $\mathcal{P}_s$ of
progressive measurable sets after time $s$ as the $\sigma$-field of
sets $P\in\mathcal{B}([s,T])\otimes\mathcal{F}_{s,T}$, the product
$\sigma$-field, such that $\chi_{P}(t,\omega)$, $t\geq s$, is a
$\mathcal{F}_{s,t}$ progressive measurable (in $t$) process, where
$\chi$ is the indicator function. We denote by $\mu_s$ the measure
 on
$\mathcal{P}_s$ defined by $\mu_s(P)=\mathbf{E}\int_s^T\chi_P(s,\omega)\,dt$.

 For definitions of consistent
processes see \citet{Londono04}. 
However for the sake of completeness we review these definitions here.
Let
$\varphi(s,t,x,\omega)$, $0\leq s \leq t\leq T$, $x\in\mathbb{D}$ be a
$\mathbb{R}^n$-valued random field on the probability space
$(\Omega,\mathcal{F},\mathbf{P})$.  We  call it a 
\emph{continuous
  $C^{m,\delta}(\mathbb{D}\colon\mathbb{R}^n)$-semimartingale} if  $\varphi_s\colon t\to\varphi(s,t,\cdot)$,  is a measurable
random field with values in
$C^{m,\delta}(\mathbb{D}\colon\mathbb{R}^n)$, that is a continuous
$(\mathcal{F}_{s,\cdot})$ semimartingale process
  decomposed
as $\varphi(s,t,x)=\varphi_{loc}(s,t,x)+\varphi_{fv}(s,t,x)$, where
$\varphi_{loc}(s,\cdot,\cdot)$ is a continuous
$C^{m,\delta}(\mathbb{D}\colon\mathbb{R}^n)$-local-martingale, and $\varphi_{fv}(s,\cdot,\cdot)$ is a continuous
$C^{m,\delta}(\mathbb{D}\colon\mathbb{R}^n)$-process of bounded
variation for each $0\leq s\leq T$.
  A pair $(a,b)$ where
$a(s,t,x,y)$ and $b(s,t,x)$ are measurable random
fields $\mathcal{F}_{s,t}$-progressive measurable in $t$, for all
$x,y\in\mathbb{D}$, $0\leq s\leq T$, is said to be the \emph{local
  characteristics of $\varphi$}, if $(a(s,\cdot,x,y),b(s,\cdot,x))$
is the local characteristic of
$\varphi_s\equiv\varphi(s,\cdot,\cdot)$ (see \citet{hK90}) for any
$s\leq t\leq T$.  In addition, a pair $(\sigma,b)$ where
$\sigma(s,t,x)$ is a measurable random field with values in
$L(\mathbb{R}^d\colon \mathbb{R}^n)$, where
$L(\mathbb{R}^d\colon\mathbb{R}^n)$ denotes the set of matrices with
size $n\times d$, $(\mathcal{F}_{s,t})$-progressive measurable in $t$,
for all $x\in\mathbb{D}$, $0\leq s\leq T$ , and $b$ is as above is
said to be the \emph{volatility and drift processes of $\varphi$} if
  \begin{equation*}
      \varphi_{loc}(s,t,x)(\omega)=\int_s^t\sigma(s,u,x)\,dW_s(u), 
  \end{equation*}
  for all $x$, $s$, $t$ and $\omega$.  If $b(s,\cdot,\cdot)$ and $\sigma(s,\cdot,\cdot)$ are
  processes of class $C^{m,\delta}$ for all $0\leq s\leq T$ , we shall say that $\varphi$ has
  \emph{volatility and drift of class $C^{m,\delta}$}.

 Let $\varphi(s,t,x)$ and $\psi(s,t,x)$ be continuous
 $C(\mathbb{D}\colon\mathbb{R}^n)$ and
 $C(\mathbb{D}\colon\mathbb{D})$ semimartingales,
  respectively; in addition, it is assumed that $\psi(s,s,x)=x$ for
  all $x\in\mathbb{D}$, and $0\leq s\leq T$.  We  say that the
  process $\varphi$ is a \emph{$\psi$-consistent semimartingale process} if for each
  $0\leq s\leq s^{\prime}\leq T$ there exists a set
  $N_{s,s^{\prime}}\in \mathcal{P}_{s^{\prime}}$ with
  $\mu_{s^{\prime}}(N_{s,s^{\prime}})=0$, such that
  $\varphi(s,t,x)=\varphi(s^{\prime},t,\psi(s,s^{\prime},x))$ for
  all $(t,\omega)\notin N_{s,s^{\prime}}$ and all $x\in\mathbb{D}$.
  We  say that the process $\varphi$ is a \emph{consistent
    semimartingale process} if $\varphi$ is a $\varphi$-consistent process.

We assume $n+1$ stocks whose \emph{evolution
  price process} $P$ is a consistent $C(\mathbb{R}_+^{n+1}\colon
\mathbb{R}_+^{n+1})$-semimartingale. 
  For $0\leq i\leq n$ we define \emph{the price per-share
  process for the $i$-stock}, $P_i$, to be the $P$-consistent
$C(\mathbb{R}_+^{n+1}\colon \mathbb{R}_+)$-semimartingale process
$P_i=\left\{P_i(s,t,p)=\pi_i\circ P(s,t,p), p\in \mathbb{R}^{n+1},
  0\leq s\leq t\leq T\right\}$ where $\pi_i$ denotes the projection
on the $i$-component.  We assume   consistent
progressive measurable $P$-consistent processes    $\sigma_{i,j}$,  $b_i$,
$\delta_i$, $r$, and $\theta_i$  of class
$C^{\delta}(\mathbb{R}_+^{n+1}\colon\mathbb{R})$ for some
$\delta>0$, where $\mathbb{R}_+$ denotes the set of real positive
numbers.    It is assumed that   $\sigma_{i,j}$,  $b_i$,
$\delta_i$, $r$, and $\theta_i$  relate to $P$ through the following stochastic
differential equations
\begin{eqnarray*}
  dP_i(s,t,p)=P_i(s,t,p)\left[b_i(s,t,p)dt + \sum_{1\leq
        j\leq d}\sigma_{ij}(s,t,p)\,dW^j_s(t)\right]  \\\nonumber
  P_i(s,s,p)=p_{i}, i=1,\ldots,n
\end{eqnarray*}
where $W^j_s(t)=W^j(t)-W^j(s)$, and,
\begin{eqnarray*}
  dP_0(s,t,p)=P_0(s,t,p)\left[-r(s,t,p)dt - \sum_{1\leq
        j\leq d}\theta_{j}(s,t,p)\,dW^j_s(t)\right]  \\\nonumber
  P_0(s,s,p)=p_0.
\end{eqnarray*}

  Throughout  this paper we shall assume that $\theta(s,\cdot,p)\in
  \ker^{\perp}(\sigma(s,\cdot,p))$,  (where $\ker^{\perp}(\sigma(s,\cdot,p)$ denotes the orthogonal
complement of the kernel of $\sigma(s,\cdot,p)$) and,
\begin{equation*}
  b(s,t,p)+\delta(s,t,p)-r(s,t,p)\mathbf{1}_n
  =\sigma(s,t,p)\theta(s,t,p)
\end{equation*}
a.e. $\mu_s$, for all $p\in\mathbb{R}_+^{n+1}$, and $0\leq s\leq T$,
where $\mathbf{1}_n^{\prime}=(1,\cdots,1)\in\mathbb{R}^n$.  This
latter assumption implies that there are not state-tame arbitrage
opportunities (see \citet{Londono04}).
  
The process of bounded variation
$B=\left\{B(s,t,p)\right\}$,
whose evolution $B(s,\cdot,p)$, $p\in\mathbb{R}_+^{n+1}$, $0\leq
s\leq T$ is given by the stochastic differential equation
\begin{equation*}
  dB(s,t,p)=B(s,t,p)r(s,t,p)dt,\qquad B(s,s,p)=1, \text{ for
  }0\leq s\leq t\leq T
\end{equation*}
will be called the \emph{bond price process}.

We shall say that $\mathcal{M}=(P,b,\sigma,\delta,r,p^0)$ is a
\emph{financial market with terminal time $T$ and initial time $0$}, 
if $b=(b_1,\ldots,b_n)$ is a vector of rate of return processes,
$\sigma=(\sigma_{i,j})$ is a matrix of volatility coefficient
processes, $\delta=(\delta_1,\ldots,\delta_n)$ is a vector of dividend
rate processes, $r$ is an interest rate process,   and 
$\theta^{\prime}(s,t,p)=\left(\theta_1(s,t),\cdots,\theta_d(s,t)\right)$
is the \emph{market price of risk}, 
and $p^0\in\mathbb{R}_+^{n+1}$ is a vector of initial prices.

We define the \emph{state price density process} to be the continuous
$C(\mathbb{R}^{n+1}_+\colon\mathbb{R}_+)$-semimartingale process
defined by
\begin{equation*}
  H(s,t,p)=B^{-1}(s,t,p)Z(s,t,p)\qquad\text{for } p\in
  \mathbb{R}_+^{n+1}, 0\leq s\leq t \leq T
\end{equation*}
where
\begin{equation*}
  Z(s,t,p)=\exp\left\{-\int_{s}^t\theta^{\prime}(s,u,p)\,dW_{s}(u)\ 
-\frac{1}{2}\int_{s}^t\left\|\theta(s,u,p)\right\|^2\,du\right\}
\end{equation*}
for $0\leq s\leq t\leq T$, and $B^{-1}(s,t,p)=1/B(s,t,p)$. 

  Assume that $\tau=\{\tau_s(x,p);
0\leq s\leq T, x\in\mathbb{R}, p\in\mathbb{R}_+^{n+1}\}$ is a
measurable family of stopping times.  
\emph{A wealth
  structure} is a triple $(X, \tau, x_0)$, where
$x_0\in\mathbb{R}$, and $X=\{X(s,t,x,p); x\in\mathbb{R},
p\in\mathbb{R}_+^{n+1}, 0\leq s\leq t\leq \tau_s(x,p)\}$, is a family of
continuous semimartingale processes with the property that
$((X,P),\tau)$ is a consistent stopping structure where $(X,P)$ is defined as
\[
(X,P)=\left\{(X(s,t,x,p),P(s,t,p)), x\in\mathbb{R}, p\in
  \mathbb{R}_+^{n+1}, s\leq t\leq \tau_s(x,p)\right\}.
\]
It is also assumed that the  drift and
volatility of the process $((X,P),\tau)$ is of  class $C^{\delta}$
for some $\delta=\delta^{X}>0$. We  say that $x_0$ is the \emph{initial value for the wealth
  process}, and we  say that $(X,\tau)$ is a \emph{wealth
  evolution structure}; we shall denote this by writing
$(X,\tau)\in\mathcal{X}(\mathcal{M})$.  
For the definition of consistent stopping structure, measurable family
of stopping times, and related ones, see
\citet{Londono2005a}. 

For the sake of completeness we also review  the definition for consistent stopping structure.
  Let $\tau=\left\{\tau(s,x), x\in\mathbb{D},0\leq s\leq T\right\}$
  be a family of stopping times with values in $\left[0,T\right]$, and
  semimartingales processes with drift and volatility of class $C^{m,\delta}$.
  It is assumed that for each $0\leq s\leq T$, $x\in\mathbb{D}$,
  $\tau(s,x)$ is a stopping time relative to the filtration
  $\left\{\mathcal{F}_{s,t}; s\leq t\leq T\right\}$ with values in $[s,T]$, and that
  $\tau(s,x)(\omega)$ is a measurable random field that  is lower
  semi-continuous with respect to $(s,x)$.  Assume that  $\psi(s,t,x), 0\leq s\leq t\leq
  \tau(s,x), x\in\mathbb{D}$ is a a family of process with values in
  $\mathbb{R}^k$, with   $\psi_{\tau}=\{\psi(s,\tau(s,x)\wedge t,x), 0\leq t\leq T ; x,
  0\leq s\leq T\}$ where $s\wedge t=\min\{s,t\}$ is a
  $C(\mathbb{D}\colon\mathbb{R}^k)$-semimartingale.   We shall say
  that $(\psi,\tau)$ is a \emph{consistent stopping structure}  if for each $0\leq s\leq T$ there exist
  $N_{s}\in\mathcal{P}_s$, $\mu_s(N_s)=0$ with
  $\tau_s(x)=\tau_{t\wedge\tau}(\psi(s,t\wedge\tau,x))$ for all
  $(t,\omega)\notin N_s$ and all $x$.  We
  say that the consistent stopping structure $(\psi,\tau)$ is of
  class $C^{m,\delta}(\mathbb{D}\colon\mathbb{R}^k)$ if $\psi_{\tau}$
  is a process of class $C^{m,\delta}(\mathbb{D}\colon\mathbb{R}^k)$.
  Given a consistent stopping structure $(\psi,\tau)$, we  say
  that a family of $\mathbb{R}^n$-valued processes
  $\varphi=\{\varphi(s,t,x), s\leq t\leq \tau_s(x) ; x\in\mathbb{D},
  0\leq s\leq T\}$ is a \emph{$\psi$-consistent process with random
    time $\tau$}, if $\varphi_{\tau}$ is a $\psi_{\tau}$-consistent
  measurable process with two parameters. Similarly,
  we  say that $\varphi$ is a process of class
  $C^{m^{\prime},\delta^{\prime}}(\mathbb{D}\colon\mathbb{R}^n)$ if
  $\varphi_{\tau}$ is a process of the same class.

 Let $\Gamma$ be a continuous semimartingale process with random time
$\tau$  with drift and volatility of class 
$C^{\delta}$
(where the positive number $\delta$ depends on $\Gamma$) with the
property that $\Gamma(s,s,x,p)=0$  and
\[
\Gamma(s,t^{\prime},x,p)+\Gamma(t^{\prime},t,X(s,t^{\prime},x,p),P(s,t^{\prime},p))=\Gamma(s,t,x,p)
\]
for all $x\in\mathbb{R}$,
$p\in\mathbb{R}_+^{n+1}$, and $0\leq s\leq \tau_s(x,p)$. 
We  say that a process $\Gamma$ as above is an \emph{income
  evolution structure for the wealth evolution structure $(X,\tau)$},
and we say that $(X,\Gamma,\tau)$ is a \emph{wealth and income
  evolution structure}.  If $\Gamma(s,t,x,p)\leq 0$ for all $x$, $p$,
$s\leq t\leq \tau_s(x,p)$ we  say that $\Gamma$ is a
\emph{consumption evolution structure for the wealth evolution
  structure $(X,\tau)$}.  Let
$(\pi_0,\pi)=\{(\pi_0(s,t,x,p),\pi(s,t,x,p)); x\in\mathbb{R},
p\in\mathbb{R}_+^{n+1}, 0\leq s\leq t\leq \tau_s(x,p)\}$ be a
$(X,P)$-consistent progressive measurable process of class
$C^{\delta}$ for some $\delta>0$ with random time $\tau$, and
$\pi_o+\pi^{\prime}\mathbf{1}_n=X$ satisfying
\begin{multline*}\label{eq:wealthportfolio_process} 
  B^{-1}(s,t,p)X(s,t,x,p)=x+\int_{s}^tB^{-1}(s,u,p)\,d\Gamma(s,u,x,p)\\
  +\int_{s}^tB^{-1}(s,u,p)\pi^{\prime}(s,u,x,p)\sigma(s,u,p)\,dW_{s}(u) \\
  +\int_s^tB^{-1}(s,u,p)\pi^{\prime}(s,u,x,p)(b(s,u,p)+\delta(s,u,p)-r(s,u,p)\mathbf{1}_n)\,du
\end{multline*}
for all $x\in\mathbb{R}$, $0\leq s \leq t\leq \tau_s(x,p)$,
$p\in\mathbb{R}^{n+1}_+$.  We  say that
$((\pi_0,\pi),\Gamma,\tau)$ as above is a \emph{portfolio evolution
  structure with random time $\tau$, financed by the
  income $\Gamma$}.   We  say that a \emph{wealth
  evolution structure $(X,\tau)\in\mathcal{X}(\mathcal{M})$  is financed by the income structure $\Gamma$}, if there exists
a portfolio evolution structure $((\pi_0,\pi),\Gamma,\tau)$ with
random time $\tau$  with
$\pi_0+\pi^{\prime}\mathbf{1}_n=X$.  In this case we  say that
$(X,\Gamma,\tau)$ is a \emph{hedgeable wealth-income structure}. 

Next we discuss the concept of utility that we shall use in this paper.
 \begin{definition}\label{d:utility}
   Consider a function $U\colon\left(0,\infty\right)\mapsto\mathbb{R}$
   continuous, strictly increasing, strictly concave and continuous
   differentiable, with
   $U^{\prime}(\infty)=\lim_{x\to\infty}U^{\prime}(x)=0$ and
   $U^{\prime}(0+)\triangleq \lim_{x \downarrow
     0}U^{\prime}(x)=\infty$.  Such a function will be called a
   utility function.
\end{definition}

Classic examples of utility functions are $U_{\alpha}(x)=
x^{\alpha}/\alpha$ for some $\alpha\in (0,1)$, $0\leq x < \infty$, and
$U(x)=\log(x)$.
For every utility function $U(\cdot)$, we shall denote by $I(\cdot)$
the inverse of the derivative $U^{\prime}(\cdot)$; both of these
functions are continuous, strictly decreasing and map $(0,\infty)$
onto itself with $I(0+)=U^{\prime}(0+)=\lim_{x\to 0^+}U^{\prime}(x)=\infty$,
$I(\infty)=\lim_{x\to\infty}I(x)=U^{\prime}(\infty)=0$. We extend $U$  by  $U(0)=U(0^+)$, and we
keep the same notation to the extension to $[0,\infty)$ of $U$  hopping that
it would be clear to the reader to which function  we are referring.  It is a well known result that 
\begin{equation}
  \label{eq:duality}
\max_{0<x<\infty}\left(U(x)-xy\right)=U(I(y))-yI(y),\qquad
  0<y<\infty  
\end{equation}

\begin{definition}\label{d:state_utility}
  Consider a  continuous function $U_1\colon
  [0,T]\times\left(0,\infty\right)\mapsto\mathbb{R}$,  such that $U_1(t,\cdot)$ is a utility function in
  the sense of Definition \ref{d:utility} for all $t\in [0,T]$.
  It follows that 
 $I_1(t,x)\triangleq(\partial
   U_1(t,x)/\partial x)^{-1}$, the inverse of the derivative of $U$,
   is a continuous  function.
 Similarly if a   utility
  function   $ U_2\colon\left(0,\infty\right)\mapsto\mathbb{R}$  is
  given then   $I_2(t,x)\triangleq (\partial
  U_2(t,x)/\partial x)^{-1}$ is continuous.  Let us denote
\begin{equation}\label{eq:inti}
\mathcal{X}(t,y)\triangleq I_2(y)+\int_t^TI_1(t^{\prime},y)\,dt^{\prime}.
\end{equation}
   We shall
  call a couple of
  functions as above   a \emph{ state preference structure}.
\end{definition}

Under the conditions outlined in the previous  definition, it is easy
to see that
$\mathcal{X}\colon [0,T]\times(0,\infty)\to(0,\infty)$ is a  continuous
function  with the property that for
each $t$, $\mathcal{X}(t,\cdot)$  maps $(0,\infty)$ onto itself,  is strictly decreasing with
$\mathcal{X}(t,0+)=\lim_{y\downarrow 0}\mathcal{X}(t,y)=\infty$ and
$\mathcal{X}(t,\infty)=\lim_{y\to\infty}\mathcal{X}(t,y)=0$.

 We extend
  $U_1$ and $U_2$ by defining 
   $U_1(t,0)=U(t,0^+)$, for all $0\leq t \leq T$ and
  $U_2(0)=U_2(0^+)$, and we keep the same notation to the extension of
  $U_1$ to $[0,T]\times [0,\infty)$, and the extension of $U_2$ to
  $[0,\infty)$.  We hope that it would be clear to the reader to
  which function we are referring.

  We point out that 
$\mathcal{X}^{-1}$ defined for each $t$ as $\mathcal{X}^{-1}(t,\cdot)$, the inverse of $\mathcal{X}(t,\cdot)$,  share the same   properties to
the some  mentioned for $\mathcal{X}$. 
We next discus the meaning of those utility functions defined above.
We should interpret  $U_1(t,x)$, for $t\in [0,T]$ the level of
``happiness'' for an agent consuming $x$ units of wealth per unit of
time at time $t$, as valued at time $0$, when the agent is planning
its consumption. Similarly, we should understand for $U_2(x)$ the
level of ``happiness'' for an agent having a final wealth of $x$ units
(at time $T$) as valued at time $0$.   This is
contrary with the traditional approach where an agent has preferences
on their consumption behavior according to their value as discounted by
a bank account, and is closer in approach to a utility function that
is state dependent.  See the literature on state dependent utilities
cited above.

For $s\leq t$  define
$\alpha(s,t)=\mathcal{X}(s,\mathcal{X}^{-1}(t,\cdot))$. Then
$\alpha(s,t)=\alpha(s,t^{\prime})\circ\alpha(t^{\prime},t)$ for all
$s$, $t$, and $t^{\prime}$ in $[0,T]$, where $\circ$ denotes standard
composition of functions.   We also observe that if $\alpha^I(s,t)\triangleq
I_1(s,\mathcal{X}^{-1}(t,\cdot))$ then $\alpha^I(s,t)\circ\alpha(t,s)=\alpha^I(s,s)$.  Throughout this paper we shall
assume the following condition on the utility structure.
\begin{condition}[Homogenity]\label{cond:homogenety}
  Let $(U_1,U_2)$ be a state preference structure defined as above.
  For all $s,t\in [0,T]$ there exist constants $\alpha_{s,t}$ and $\alpha_s^I$ such that
  $\alpha(s,t)(x)=\alpha_{s,t}x$, and $\alpha^I(s,s)(x)=\alpha^I_sx$
  where $\alpha(s,t)$ and $\alpha^I(s,t)$ are defined as
  the previous paragraph.  In this case we say that $(U_1,U_2)$ is a
  \emph{homogeneous state preference structure}.
\end{condition}
A way to see this  is to say that  the
structure for the utility preferences remains the same  as
time evolves.  We next describe some important
examples that fit the previous conditions. 
\begin{example}\label{ex:xtoalpha}
  Assume a continuous positive  function
  $h\colon[0,T]\to(0,\infty)$, and assume that
  $U_1(t,x)=x^{\alpha}h(t)$ and $U_2(x)=cx^{\alpha}$ with
  $\alpha\in(0,1)$ and $c\geq 0$.  This is an state preference structure that
  satisfies Condition~\ref{cond:homogenety}. Indeed in this case 
\[
\alpha_{s,t}=\frac{c^{1/(1-\alpha)}+\int_s^Th^{1/(1-\alpha)}(t^{\prime})\,dt^{\prime}}{c^{1/(1-\alpha)}+\int_t^Th^{1/(1-\alpha)}(t^{\prime})\,dt^{\prime}},\qquad
\alpha_t^I=\frac{h^{1/(1-\alpha)}(t)}{c^{1/(1-\alpha})+\int_t^Th^{1/(1-\alpha)}(t^{\prime})\,dt^{\prime}}
\]
\end{example}
\begin{example}\label{ex:log}
Assume a continuous positive function
  $h$ as above, and assume that
  $U_1(t,x)=h(t)\log(x)$ and $U_2(x)=c\log(x)$, with $c\geq 0$.  
It follows that this is a state preference structure that satisfies
Condition~\ref{cond:homogenety} with
\[
\alpha_{s,t}=\frac{c+\int_s^Th(t^{\prime})\,dt^{\prime}}{c+\int_t^Th(t^{\prime})\,dt^{\prime}},\qquad
\alpha_t^{I}=\frac{h(t)}{c+\int_t^Th(t^{\prime})\,dt^{\prime}}
\]
\end{example}
\begin{example}\label{ex:separability}
  Let $U_1(t,x)=h(t)u(x/h(t))$  and $U_2(x)=cu(x/c)$, where $u(\cdot)$
  is a utility function, $h(\cdot)$ is a positive continuous function
  and $c>0$.
 It follows that   $(U_1,U_2)$ is a state preference structure that
 satisfies Condition~\ref{cond:homogenety}.  In this case 
\[
\alpha_{s,t}=\frac{c+\int_s^Th(t^{\prime})\,dt^{\prime}}{c+\int_t^Th(t^{\prime})\,dt^{\prime}},\qquad
\alpha_t^I=\frac{h(t)}{c+\int_t^Th(t^{\prime})\,dt^{\prime}}
\]
In particular, when $h\equiv 1$ we obtain that
$U_1(t,x)=u(x)$, and $U_2(x)=cu(c^{-1}x)$ for some $c>0$ define a
state preference structure that satisfies Condition~\ref{cond:homogenety}.
\end{example}
\section[A Martingale approach]{A Martingale approach}
\label{sec:martingale-approach}

\begin{definition}\label{d:ConsumptionWealth}
 Assume that
  $(X,\Gamma,\tau)\in\mathcal{X}(\mathcal{M})$ is a wealth-income
  evolution structure.   Assume that
  $\Gamma\equiv E-C$, where
  \begin{equation}
    \label{eq:cumulative_consumption}
dC(s,t,x,p)=c(s,t,x,p)\,dt    
  \end{equation}
  and,
\begin{equation}
  \label{eq:cumulative_endowment}
dE(s,t,p)=\varepsilon(s,t,p)\, dt   
\end{equation}
for non-negative $(X,P)$ consistent processes $c$ and $\varepsilon$ of
class $C^{\delta}$ for some $\delta>0$.
Moreover assume that \[
\mathbf{E}\left[\int_{s}^TH(s,u,p)\varepsilon(s,u,p)\, du\right]<\infty
\]
for all $p\in\mathbb{R}_+^{n}$, and $0\leq s\leq T$.  Similarly
assume that 
\[
\mathbf{E}\left[\int_{s}^TH(s,u,p)c(s,u,x,p)\, du\right]<\infty
\] 
for all $p\in\mathbb{R}_+^{n}$, $x\in\mathbb{R}$ and $0\leq s\leq T$.

We should say $(X,c,\varepsilon,\tau)$ as above is a \emph{rate of
  consumption and endowment evolution structure}.  We shall say that
$c$ is the \emph{consumption rate evolution structure}, and
$\varepsilon$ is the \emph{endowment rate evolution structure}.  We
also say that $E$ is a \emph{cumulative endowment structure}, and $C$ is a
\emph{cumulative consumption structure}.

A \emph{minimal wealth structure} $L$ is a $P$ consistent process with
drift and volatility of class
$C^{\delta}$ for some $\delta>0$ where
$L(s,\cdot,p)H(s,\cdot,p)$ is uniformly bounded below for all $p$, $s$,
(where the bound might depend on $p$ and $s$) 
such that 
\[
\mathbf{E}\left[H(s,t,p)L(s,t,p)\right]<\infty
\]
for all $p$, $s$ and $t$.
\end{definition}

We should emphasize that the name might be confusing, since a minimal
wealth structure \emph{is not} a wealth structure in the sense defined
above.  However we keep the name for consistency  with standard use.  It is natural to believe that the evolution of income due to labor
only depends on the evolution of the state of the economy and not on the
current wealth of an agent.   

Typically we are interested  in consumption and endowment structure
evolution structures whose wealth remains above some given process.
Next we present the definition that embodies this idea.

\begin{definition}\label{def:portfolio_consumption_admisible}
  Let $(X,\varepsilon,c,T)$ (denoted by $(X,\varepsilon,c)$) be a
  hedgeable (by a state tame portfolio) cumulative consumption and endowment
  evolution structure,   as in Definition
  \ref{d:ConsumptionWealth}, with portfolio evolution structure
  $(\pi_0,\pi)$.  We shall say that the couple $(\pi,c)$ of portfolio
  on stocks and rate of consumption, is
  \emph{admissible for $(L,\varepsilon)$}  and write
  $(\pi,c)\in\mathcal{A}(L,\varepsilon)$ if for any   $x$, $s$
  and $p$ with
  $x\geq L(s,s,p)$
\begin{equation}\label{eq:goal}
  X(s,t,x,p)\geq L(s,t,p)\qquad\mbox{a.e. }  
\end{equation}  
If there is not portfolio
  on stocks and rate of consumption for
$(L,\varepsilon)$  we should say
that the class cited  above is empty, and we would denote this by
$\mathcal{A}(L,\varepsilon)=\emptyset$
\end{definition}
 For any hedgeable
  wealth and income evolution structure $(X,E-C)$ with $(\pi,c)$
  admissible for  $(L,\varepsilon)$  it must hold that 
\[
x\geq
\mathbf{E}\left[H(s,T,p)L(s,T,p)+\int_{s}^T{H(s,u,p)\left(c(s,u,x,p)-\varepsilon(s,u,p)\right)}{du}\right]
\]
for any $x\geq L(s,s,p)$, 
where the latter follows since  the process defined by
equation (\ref{eq:localmartingale}) is a super-martingale.  It is often
the case that $L(s,T,p)=0$ for all $s$ and $p$.  In these latter case
the condition for the previous equation becomes
\begin{equation*}
x\geq
\mathbf{E}\left[\int_{s}^T{H(s,u,p)\left(c(s,u,x,p)-\varepsilon(s,u,p)\right)}{du}\right].
\end{equation*}

Next we explain the problem that we are interested to solve in this
paper.  We assume  a minimal wealth structure
$L$ and an endowment rate evolution structure $\varepsilon$.  The
control stochastic problem that we propose to solve concerns a small
investor that at time $0$ has an initial capital $x$,  is
constrained to not  let his wealth to fall below a minimal wealth
process $L(0,\cdot,p)$, has a rate of
endowment process, $\varepsilon(0,\cdot,p)$ and has
at his disposal portfolio/consumption processes
$(\pi,c)\in\mathcal{A}(L,\varepsilon)$. The following
Proposition \ref{thm:existence_portfolios}   is a direct consequence of \citet[Theorem 2]{Londono2005a}; it
provides conditions under which $\mathcal{A}(L,\varepsilon)\neq\emptyset$.

\begin{proposition}\label{thm:existence_portfolios}
  Assume a  minimal wealth structure $L$, and a
  rate evolution structure  $\varepsilon$ as in Definition
  \ref{def:portfolio_consumption_admisible}.  Assume that 
\begin{equation*}
H(s,t,p)L(s,t,p)- \int_s^tH(s,t,p)\varepsilon(s,t,p)\,du
\end{equation*}
is a martingale for all $s$, $p$. 
Then, there exist a cumulative consumption and endowment evolution structure
$(X,0,\varepsilon)$  with
$(\pi,0)\in\mathcal{A}(L,\varepsilon)$ where $\pi$ is the
portfolio on stocks defined by $X$.
\end{proposition}
\begin{proof}
Define $X$ by
\[
X(s,t,x,p)\triangleq
 L(s,t,p)+  (x-L(s,s,p))H^{-1}(s,t,p)    
\]
It follows using \citet[Theorem 2]{Londono2005a} that
$(X,0,\varepsilon)$ is the desired cumulative consumption and
endowment evolution structure. 
\end{proof}

The following condition is needed in order to solve the problem of
optimal investment and consumption under less strict conditions on the
minimal wealth structure. For the following condition let  $(X,c,\varepsilon,\tau)$ be  a rate of consumption and endowment evolution structure
with \emph{discounted payoff process} defined as 
   \begin{equation}
      \label{eq:localmartingale}
Y(s,t,x,p)\triangleq H(s,t,p)X(s,t,x,p) +
\int_{s}^t{H(s,u,p)\left(c(s,u,x,p)-\varepsilon(s,u,p)\right)}{du}.      
    \end{equation}
  \begin{condition}\label{con:continuity}
Let  $(X,c,\varepsilon,\tau)$ be  a rate of consumption and endowment
evolution structure, as above. 
We assume that  for all stopping times
    $\tau\in\mathcal{S}(X)$, and $0\leq s\leq T$ the function
\begin{equation*}
\varphi_{s,\tau}(x,p)=\mathbf{E}\left[Y(s,t\wedge\tau(s,x,p),x,p)\right] 
\end{equation*}
is a continuous function in $(x,p)$, and the given family of functions
is an equicontinuous set of functions on compact sets (in $(x,p)$),
where $s\wedge t=\min(s,t)$ and it is assumed that
$\sup\emptyset =\infty$.
Here $\mathcal{S}(X)$ denotes,  the family of stopping times
that are  $(X,P)$-consistent.
Moreover assume that there exist positive constants $\gamma\geq 1$,
$\alpha_1,\alpha_2,\alpha_3,\beta_0,\cdots,\beta_n$, with
$\alpha_1^{-1}+\alpha_2^{-1}+\alpha_3^{-1}+\sum_{i=0}^n\beta_i<1$ such
that the random field $Y(s,t,x,p)$
satisfies
\begin{multline*}
  \mathbf{E}\left[\mid
    Y(x,p,s,t)-Y(x^{\prime},p^{\prime},s^{\prime},t^{\prime})\mid^{\gamma} \right]\leq \\
  C\left(\mid s-s^{\prime}\mid^{\alpha_1}+\mid
    t-t^{\prime}\mid^{\alpha_2}+ \mid x-x^{\prime}\mid^{\alpha_3}
    +\sum_{i=0}^n\mid p_i-p_i^{\prime}\mid^{\beta_i}\right).
\end{multline*}
  \end{condition}
This condition is usually satisfied when $X$ is a process
  that solves a stochastic differential equation.  For instance, see
  \citet[Lemma
  4.5.6]{hK90}.  The  last inequality  is needed
  in order to obtain a continuous modification of the random field and
  its conditional expectation.  See Kolmogorov's continuity criterion
  for random fields (\citet[Theorem 1.4.1 and Exercise 1.4.12]{hK90}).
  In the following, conditional expectations of stochastic processes
  are the continuous modifications of the given stochastic processes. 

For the problems of optimal consumption and terminal wealth that we
describe below we shall assume that the minimal wealth structure is
defined in a way such that the discounted minimal wealth process of an
agent can not fall below the current value of future endowments,
\begin{equation}
  \label{eq:future_endownment}
L(s,t,p)=\frac{-1}{H(s,t,p)}\mathbf{E}\left[\int_t^TH(s,u,p)\varepsilon(s,u,p)\,du\mid\mathcal{F}_{s,t}\right]  ,
\end{equation}
for $0\leq s\leq t\leq T$,  $p\in\mathbb{R}_+^{n}$.
In fact is not difficult to see that the family of stochastic
processes defined by the last equation, is a minimal wealth
process with $\mathcal{A}(L,\varepsilon)\neq\emptyset$,  since the discounted payoff process
$Y(s,t,p)$ satisfies
\begin{multline*}
  Y(s,t,p)= H(s,t,p)L(s,t,p)-\int_s^tH(s,u,p)\varepsilon(s,u,p)\\
  = -\mathbf{E}\left[\int_{s}^{T}H(s,u,p)\varepsilon(s,u,p)\,du
    \mid\mathcal{F}_{s,t}\right]
\end{multline*}
and therefore is clearly a martingale. In fact,
Proposition~\ref{thm:existence_portfolios} is a consequence of a more
general theorem stated  below. It  allows to solve the problem of
optimal consumption and investment under more general minimal wealth
structures.

  \begin{theorem}\label{thm:admisibility}
    Let $(X,c,\varepsilon,\tau)$ be a rate of consumption and endowment evolution
    structure as in
    Definition~\ref{d:ConsumptionWealth} with cumulative endowment and
    consumption structures $C$, and $E$ as defined  by equations
    ~\eqref{eq:cumulative_consumption} and
    ~\eqref{eq:cumulative_endowment} respectively. If the family of
    processes defined by equation \eqref{eq:localmartingale}
are martingales for each $x$, $p$,
$s$ then $(X,E-C)$ is a hedgeable wealth-income structure. Moreover,  if
 $Y(s,\cdot,x,p)$ is a super-martingales for each 
$x$, $p$ and $s$,  and  Condition~\ref{con:continuity} holds,  then
$(X,E-C)$ is dominated by a hedgeable wealth-income structure
$(X^{\prime},E-C)$.  See \citet{Londono2005a} for the definition of
this concept.
In particular, if  a  minimal wealth evolution structure $L$  is given
such that 
\[
H(s,t,p)L(s,t,p) -
\int_{s}^t{H(s,u,p)\varepsilon(s,u,p)}{du}      
\]
is a super-martingale for all $s$,  and $p$ with the property that for
all $x\in\mathbb{R}$
\begin{equation}\label{eq:initial_minimal_wealth}
x\geq L(s,s,p)+\mathbf{E}\left[\int_s^T{H(s,u,p)c(s,u,x,p)}\,{du}\right],
\end{equation}
then there exist a hedgeable wealth and income evolution structure
$(X^{\prime},E-C)$ with portfolio on stocks $\pi^{\prime}$, such that
$X^{\prime}(s,t,x,p)\geq L(s,t,p)$ for all $s$, $p$, $t$, and $x$, 
satisfying  equation (\ref{eq:initial_minimal_wealth}).  In  particular  $(\pi^{\prime},c)\in\mathcal{A}(L,\varepsilon)$.
  \end{theorem}
  \begin{proof}
  The first part of the proof is a straightforward consequence of
  \citet[Theorem 2 and Theorem 3]{Londono2005a}.  For the second part of
  the proof we observe that if $X$ is defined by 
\begin{multline*}
 X(s,t,x,p)=L(s,t,p)
+\frac{1}{H(s,t,p)}\mathbf{E}\left[\int_t^T{H(s,u,p)c(s,u,x,p)}{du}\mid\mathcal{F}_{s,t}\right]
+\\
\frac{1}{H(s,t,p)}\left(x-L(s,s,p)-\mathbf{E}\left[\int_s^T{H(s,u,p)c(s,u,x,p)}\,{du}\right]\right),
\end{multline*}
then $(X,E-C)$ is a wealth income evolution structure with the
property that the process defined by the equation (\ref{eq:localmartingale})
above, is a super-martingale for any $x$, $p$, and $s$.  As a
consequence of the above  there exist a hedgeable wealth
and income evolution structure $(X^{\prime},E-C)$ dominating  $(X,E-C)$.  It follows that for 
any $x$ satisfying  equation (\ref{eq:initial_minimal_wealth}),  $X^{\prime}(s,t,x,p)\geq L(s,t,p)$ for all
$p$, $s$ and $t$. 
  \end{proof}
  
It is clear from the previous theorem that it is still possible to
consider minimal wealth structures more general that the ones  we
pursue in this paper (those that satisfy equation (\ref{eq:future_endownment})).

  \section[Consumption and Portfolio Optimization]{Consumption and Portfolio Optimization}
  \label{sec:cons-portf-optim}

In this paper we are interested in solving
  the optimization problems presented in this section.  We assume a
  state preference structure $(U_1,U_2)$.  We also assume a 
  minimal wealth-income structure  $L$,  defined as in
  equation  (\ref{eq:future_endownment}), with  an  endowment rate evolution structure
  $\varepsilon$. 

\begin{problem}[Utility from consumption]\label{p:utility_consumption}
 Under the hypotheses assumed
  here \emph{the problem of maximization of utility from consumption}
  is defined to be the problem of maximizing of expected utility of
  discounted consumption,
\begin{equation*}\label{eq:value_consumption}
V_1(x,p)\triangleq\sup_{(\pi,c)\in\mathcal{A}_1(L,\varepsilon,x)}\mathbf{E}\int_{0}^T{U_1(t,H(0,t,p)c(0,t,x,p))}{\,dt} ,
\end{equation*}
for all  $p$ and
$x>-\mathbf{E}\left[\int_0^TH(0,u,p)\varepsilon(0,u,p)\,du\right]$,
where 
\[
\mathcal{A}_1(L,\varepsilon,x)\triangleq \left\{(\pi,c)\in\mathcal{A}(L,\varepsilon)\colon\mathbf{E}\int_{0}^T{U_1^{-}(t,H(0,t,p)c(0,t,x,p))}{\,dt}<\infty\right\}
\]
and $U_1^{-}(t,x)=-(U_1(t,x)\wedge 0)$.
We shall say that $V_1$ is the \emph{value function for the problem of
optimization of utility from consumption}. 
\end{problem}
\begin{problem}[Utility from terminal wealth]\label{p:utility_terminal}
   Under the hypotheses
  assumed in this section the \emph{ problem of maximization of
    utility from terminal wealth}
  is defined to be the problem of maximizing the  expected utility from
  discounted terminal wealth at time $T$,
  \begin{equation*}
  V_2(x,p) \triangleq \sup_{(\pi,c)\in\mathcal{A}_2(L,\varepsilon,x)}
  \mathbf{E}\left[ U_2(H(0,T,p)X(0,T,x,p))\right] ,    
  \end{equation*}
for all  $p$ and
$x>-\mathbf{E}\left[\int_0^TH(0,u,p)\varepsilon(0,u,p)\,du\right] $,
where 
\[
\mathcal{A}_2(L,\varepsilon,x)\triangleq \left\{(\pi,c)\in\mathcal{A}(L,\varepsilon)\colon\mathbf{E}\left[ U_2^-(H(0,T,p)X(0,T,x,p))\right]<\infty\right\}
\]
and $U_2^{-}(x)=-(U_2(x)\wedge 0)$.
We shall say that $V_2$ is the \emph{value function for the problem of
optimization of utility from terminal wealth}. 
\end{problem}
\begin{problem}[Utility from both consumption and terminal wealth]\label{p:utility_consumption_terminal}
  Under the hypotheses assumed
  above the \emph{problem of maximization of utility from both consumption and
    terminal wealth} is defined to be the problem of maximization of 
  expected utility from consumption and terminal wealth,
\begin{multline*} 
  V(x,p) \triangleq \sup_{(\pi,c,x)\in\mathcal{A}(L,\varepsilon,x)}
  \mathbf{E}[\int_{0}^{T}U_1(t,H(0,t,p)c(0,t,x,p))\, dt \\
  + U_2(H(0,T,p)X(0,T,x,p)) ]
\end{multline*}
for all  $p$ and
$x>-\mathbf{E}\left[\int_0^TH(0,u,p)\varepsilon(0,u,p)\,du\right]$,
where 
\[
\mathcal{A}(L,\varepsilon,x)\triangleq \mathcal{A}_1(L,\varepsilon,x)\cap\mathcal{A}_2(L,\varepsilon,x).
\]
We shall say that $V$ is the \emph{value function for the problem of
optimization of utility from consumption and terminal wealth}. 
\end{problem}
A few words are needed here.  With the help of Proposition
\ref{thm:admisibility} it is also possible to consider more general
optimization problems when the restriction on minimal wealth is not
necessarily the current value of future endowments as defined by
equation (\ref{eq:future_endownment}).  However in this paper we do not pursue this
line of research, and we believe this is  an open area of
research were more general restrictions on minimal wealth could be
studied.

We also point out that $0$ does not play any  special role, and the concepts like
wealth, cumulative income, portfolio process, state preference
structure, value functions  and alike can be carried out for any time interval $[s,T]$
with $0\leq s\leq T$. The above remark allows us to consider parameterized utility
preference structures with parameter  $0\leq s\leq T$, 
defined 
on the  time interval $[s,T]$.  This  models
how an agent can change preferences as time evolves. In
\citet{Londono2006} we study the investment and consumption behavior of
agents that change preferences as the result of aging.

The problems consider above are different from  the standard problems of
optimal consumption and investment, see for instance
(\citet{Karatzas98}).  First,  the optimization problems are over
portfolio and consumptions which  are \emph{consistent}.  Second, it looks at
utility functions as reflecting the level of satisfaction over levels
of  consumption
 in Problem \ref{p:utility_consumption}, final wealths in
Problem \ref{p:utility_terminal}, and on both in Problem
\ref{p:utility_consumption_terminal}, as valued by the market when
the agent is making his consumption and investment decisions (at time $0$).


Let us define
\begin{equation}\label{eq:utility_structure}
\Pi(s,t,p) \triangleq -\mathbf{E}\left[\int_t^TH(s,u,p)\varepsilon(s,u,p)\,du\mid\mathcal{F}_{s,t}\right]
\end{equation}
For any $x>\Pi(t,t,p)$ we  define $\mathcal{Y}(t,x,p)$ as the unique
solution of  
\[
\mathcal{X}(t,\mathcal{Y}(t,x,p))=x-\Pi(t,t,p)
\]
where $\mathcal{X}$ is defined by equation \eqref{eq:inti}.
 It follows that
$\mathcal{Y}(t,x,p)=\mathcal{X}^{-1}(t,x-\Pi(t,t,p))$.

\begin{theorem}\label{thm:optimal_consumption_investment}
  Assume the hypotheses of Problem
  \ref{p:utility_consumption_terminal}, and in addition   assume that $(U_1,U_2)$ is a homogeneous state preference structure
  (see Condition \ref{cond:homogenety}) .  Define $\xi$ as
\[\xi(s,t,x,p)\triangleq
\begin{cases}
  H^{-1}(s,t,p)\left(\Pi(s,t,p)+\mathcal{X}(t,\mathcal{Y}(s,x,p))\right)&\text{if $x>\Pi(s,s,p)$} \\
H^{-1}(s,t,p)\left(\Pi(s,t,p)+x-\Pi(s,s,p)\right)& \text{otherwise, }
\end{cases}
\]
and let $c$   be defined as
\[
c(s,t,x,p)\triangleq
\begin{cases}
 H^{-1}(s,t,p)I_1(t,\mathcal{Y}(s,x,p))&  \text{if $x>\Pi(s,s,p)$} \\
0 & \text{otherwise.}
\end{cases}
\]
Then, $(\xi,c,\varepsilon)$ is a hedgeable cumulative
consumption and endowment structure, with portfolio
$(\pi,c)\in\mathcal{A}(L,\varepsilon)$ that is optimal for the problem
of optimal consumption and investment.  The value function is given by 
\begin{equation*}
  V(x,p)=G(0,\mathcal{Y}(0,x,p)),
\end{equation*}
where 
\[
G(s,y)=\int_s^TU_1(t,I_1(t,y))\,dt + U_2(I_2(y))
\]
for 
$0<y <\infty$.
The corresponding optimal
portfolio on stocks is  
\begin{multline}\label{eq:optimal_portfolio}
\left(\xi(s,t,x,p)+\varPi(t,P(s,t,p))+\phi_0(t,P(s,t,p))\right)(\sigma\sigma^{\prime})^{-1}(b+\delta-r\mathbf{1}_n)(s,t,p)-\\
(\phi_1(t,P(s,t,p),\cdots,\phi_n(t,P(s,t,p)))^{\prime}
\end{multline}
where 
\[
\varPi(t,p)\triangleq\Pi(t,t,p), \qquad \phi_i(t,p)\triangleq
p_i\frac{\partial \varPi(t,p)}{\partial p_i}\qquad 0\leq i\leq n
\]
\end{theorem}
\begin{proof}
Let us point out that  Condition \ref{cond:homogenety} implies that
 $\xi$ is a (consistent)  process,
and clearly it is Lipschitz continuous.  The homogeneity also implies
that  $c$ is a $(\xi,P)$ consistent process of class $C^{0,1}$.    We observe that 
\begin{multline}\label{eq:optimal_supermartingale}
Y(s,t,x,p)\triangleq H(s,t,p)\xi(s,t,x,p)+\int_s^tH(s,u,p)(c(s,u,x,p)-\varepsilon(s,u,p))\,du\\
=x+\mathbf{E}\int_s^TH(s,u,p)\varepsilon(s,u,p)\,du-\mathbf{E}\left[\int_s^TH(s,u,p)\varepsilon(s,u,p)\,du\mid\mathcal{F}_{s,t}\right]
\end{multline}
is a martingale, and therefore Theorem \ref{thm:admisibility}
implies that $(\xi,c,\varepsilon)$ is a cumulative consumption and
endowment structure with portfolio
$(\pi,c)\in\mathcal{A}(L,\varepsilon)$.  Next, we observe that for $x>\Pi(s,s,p)$
\begin{multline*}
\mathbf{E}[\int_{s}^{T}U_1(t,H(s,t,p)c(s,t,x,p))\, dt] \\
  + \mathbf{E}[U_2(H(s,T,p)\xi(s,T,x,p))]
=G(s,\mathcal{Y}(s,x,p))  
\end{multline*}
and, if $(X^{\prime},\epsilon,c^{\prime})$ is a hedgeable rate of
consumption, endowment  and wealth evolution structure, then for $x>\Pi(s,s,p)$,
{\small
\begin{multline*}
  \mathbf{E}[\int_{s}^{T}U_1(t,H(s,t,p)c^{\prime}(s,t,x,p))\, dt+
  U_2(H(s,T,p)X^{\prime}(s,T,x,p))]\leq\\
G(s,\mathcal{Y}(s,x,p))-
\mathcal{Y}(s,x,p)\left[\int_{s}^TI_1(t,\mathcal{Y}(s,x,p))\,dt
+ I_2(\mathcal{Y}(s,x,p))\right]\\
+\mathcal{Y}(s,x,p)\mathbf{E}\left[H(s,T,p)X^{\prime}(s,T,x,p) +
\int_{s}^T{H(s,u,p)\left(c^{\prime}(s,u,x,p)\right)}{du}\right]\\
\leq G(s,\mathcal{Y}(s,x,p))
\end{multline*}}
where the first inequality is a consequence of equation
\eqref{eq:duality} and  the last inequality is a consequence of the fact that the process
defined by equation \eqref{eq:localmartingale} is a super-martingale for any hedgeable
wealth-income structure.  

Next we prove  that the optimal portfolio satisfies equation
(\ref{eq:optimal_portfolio}). It is known that the corresponding optimal portfolio should satisfy 
\[
 \sigma^{\prime}(s,t,p)\pi(s,t,x,p)=H^{-1}(s,t,p)\varphi(s,t,x,p)+ \xi(s,t,x,p)\theta(s,t,p)
\] 
where $\varphi(s,t,x,p)$ is the  process such that
\[
 Y(s,t,x,p)
  =x + \int_s^t\varphi ^{\prime}(s,u,x,p)\,dW_s(u)
\]
Using the uniqueness of the decomposition of a continuous semimartingale
as local martingale and a process of bounded variation,  Ito's rule,
and 
the fact that $\varepsilon$ is a $P$ consistent process,  it
follows by a straightforward computation that the optimal
portfolio is given by equation (\ref{eq:optimal_portfolio}).
\end{proof}
\begin{remark}\label{rem:portfolio_under_income}
If $\varepsilon$ is a $P$ consistent process where
$\mathbf{E}\left[\int_s^TH(s,u,p)\varepsilon(s,u,p)\,du\mid\mathcal{F}_{s,t}\right]$
is a deterministic function  then  the proof of the above theorem
shows that the optimal portfolio is  
\[
\pi(s,t,x,p)=(\sigma\sigma^{\prime})^{-1}(b+\delta-r\mathbf{1}_n)(s,t,p)\xi(s,t,x,p)
\]
One important example of the above case is when there are not additional
income to invest in the portfolio. 
\end{remark}

\begin{remark}\label{rem:homogenity}
 One of the consequences of Theorem
\ref{thm:optimal_consumption_investment}, is that the solution to
Problem \ref{p:utility_consumption_terminal} above, under the 
hypothesis that the state preference structure is homogeneous, is
also homogeneous in the sense that we explain next.  For any time $0\leq
s\leq T$ the solution $(\xi,c)$ of Problem
\ref{p:utility_consumption_terminal} (as well as its associated
optimal portfolio) satisfies the property that its
restriction to the time interval $[s,T]$ is also optimal for the
problem of optimal consumption and investment after time $s$ (where
the definition of the solution to the problem has been outlined after
the definition of the solution to  the problems of optimal consumption
and investment).  The latter
remark is a consequence of the proof of Theorem
\ref{thm:optimal_consumption_investment}.  
\end{remark}

Next, we present without proof the solution to the problem of optimal
consumption and investment when there is not any preference on partial
consumption of final wealth.  The proof is similar to
the proof of Theorem \ref{thm:optimal_consumption_investment}, and is 
left to the reader.  In order to state the  following theorem we introduce the functions
\[
\mathcal{X}_1(t,y)=\int_t^TI_1(t^{\prime},x)\,dt^{\prime},\qquad \mathcal{X}_2(y)=I_2(y)
\]
and
\[
G_1(t,x)=\int_s^TU_1(t,I_1(t^{\prime},y))\,dt^{\prime}\qquad G_2(y)=U_2(I_2(y))
\]
for $y>0$ and $0\leq t\leq T$.  We also set
$\mathcal{Y}_1(t,x,p)=\mathcal{X}_1^{-1}(t,x-\Pi(t,t,p))$, and
$\mathcal{Y}_2(t,x,p)=\mathcal{X}_2^{-1}(x-\Pi(t,t,p))$ for $x>0$ and
$0\leq t\leq T$.
\begin{theorem}
  Assume a homogeneous state preference structure and minimal
  wealth-income structure $(U_1,U_2)$,  $L$ defined as in  equation  (\ref{eq:future_endownment}), with  an  endowment rate evolution structure  $\varepsilon$. Define $\xi_1$ as
\[\xi_1(s,t,x,p)\triangleq
\begin{cases}
  H^{-1}(s,t,p)\left(\Pi(s,t,p)+\mathcal{X}_1(t,\mathcal{Y}_1(s,x,p))\right)&\text{if $x>\Pi(s,s,p)$} \\
H^{-1}(s,t,p)\left(\Pi(s,t,p)+x-\Pi(s,s,p)\right)& \text{otherwise,}
\end{cases}
\]
and let $c_1$ be  defined as
\[
c_1(s,t,x,p)\triangleq
\begin{cases}
 H^{-1}(s,t,p)I_1(t,\mathcal{Y}_1(s,x,p))&  \text{if $x>\Pi(s,s,p)$} \\
0 & \text{otherwise.}
\end{cases}
\]
Then, $(\xi_1,c_1,\varepsilon)$ is a cumulative
consumption and endowment structure, with portfolio
$(\pi_1,c_1)\in\mathcal{A}(L,\varepsilon)$ that is optimal for the problem
of optimal consumption and investment.  The value function is given by 
\begin{equation*}
  V_1(x,p)=G_1(0,\mathcal{Y}_1(0,x,p)),
\end{equation*}
for 
$0<y <\infty$, and optimal portfolio on stocks 
 \begin{multline*}
\left(\xi_1(s,t,x,p)+\varPi(t,P(s,t,p))+\phi_0(t,P(s,t,p))\right)(\sigma\sigma^{\prime})^{-1}(b+\delta-r\mathbf{1}_n)(s,t,p)-\\
(\phi_1(t,P(s,t,p),\cdots,\phi_n(t,P(s,t,p)))^{\prime}, 
\end{multline*}
where $\Pi$, and $\phi_i$, $0\leq i\leq n$ are defined as in Theorem \ref{thm:optimal_consumption_investment}.
\end{theorem}

\begin{theorem}
 Assume a homogeneous state preference structure and minimal wealth-income structure,  $L$ defined as in  equation  (\ref{eq:future_endownment}), with  an  endowment rate evolution structure  $\varepsilon$. Let $\xi_2$ be defined as
\[\xi_2(s,t,x,p)\triangleq
\begin{cases}
  H^{-1}(s,t,p)\left(\Pi(s,t,p)+\mathcal{X}_2(t,\mathcal{Y}_2(s,x,p))\right)&\text{if $x>\Pi(s,s,p)$} \\
H^{-1}(s,t,p)\left(\Pi(s,t,p)+x-\Pi(s,s,p)\right)& \text{otherwise.}
\end{cases}
\]
Then, $(\xi_2,0,\varepsilon)$ is a cumulative
consumption and endowment structure, with portfolio
$(\pi_2,0)\in\mathcal{A}(L,\varepsilon)$ that is optimal for the problem
of optimal consumption and investment.  The value function is given by 
\begin{equation*}
  V(x,p)=G(\mathcal{Y}_2(0,x,p)),
\end{equation*}
for 
$0<y <\infty$, and the optimal portfolio on stocks is  
 \begin{multline}
\left(\xi_2(s,t,x,p)+\varPi(t,P(s,t,p))+\phi_0(t,P(s,t,p))\right)(\sigma\sigma^{\prime})^{-1}(b+\delta-r\mathbf{1}_n)(s,t,p)-\\
(\phi_1(t,P(s,t,p),\cdots,\phi_n(t,P(s,t,p)))^{\prime}
\end{multline}
where $\Pi$, and $\phi_i$, $0\leq i\leq n$ are defined as in Theorem
\ref{thm:optimal_consumption_investment}.
\end{theorem}

\bibliographystyle{plainnat}
\bibliography{./../biblio/utility,./../biblio/eafit,./../biblio/tameness}

\begin{thebibliography}{43}
\expandafter\ifx\csname natexlab\endcsname\relax\def\natexlab#1{#1}\fi
\expandafter\ifx\csname url\endcsname\relax
  \def\url#1{{\tt #1}}\fi

\bibitem[Balduzzi and Lynch(1999)]{Balduzzi_Lynch99}
Pierlugi Balduzzi and A.~Lynch.
\newblock Transaction costs and predictability: Some utility cost calculations.
\newblock {\em Journal of Financial Economics}, 52:\penalty0 47--78, 1999.

\bibitem[Brennan(1998)]{Brennan98}
Michael~J. Brennan.
\newblock The role of learning in dynamic portfolio decisions.
\newblock {\em European Finance Review}, 1:\penalty0 295--306, 1998.
\newblock Computational financial modelling.

\bibitem[Brennan et~al.(1997)Brennan, Schwartz, and Lagnado]{Brennan_etal97}
Michael~J. Brennan, Eduardo~S. Schwartz, and Ronald Lagnado.
\newblock Strategic asset allocation.
\newblock {\em J. Econom. Dynam. Control}, 21\penalty0 (8-9):\penalty0
  1377--1403, 1997.
\newblock ISSN 0165-1889.
\newblock Computational financial modelling.

\bibitem[Brennan and Xia(2001)]{Brennan_Xia01}
Michael~J. Brennan and Yihong Xia.
\newblock Dynamic asset allocation under inflation.
\newblock {\em Journal of Finance}, 57:\penalty0 1201--1238, 2001.

\bibitem[Campbell et~al.(2004)Campbell, Chacko, Rodriguez, and
  Viceira]{Campbell_etal04}
John~Y. Campbell, George Chacko, Jorge Rodriguez, and Luis~M. Viceira.
\newblock Strategic asset allocation in a continuous-time {VAR} model.
\newblock {\em J. Econom. Dynam. Control}, 28\penalty0 (11):\penalty0
  2195--2214, 2004.
\newblock ISSN 0165-1889.

\bibitem[Campbell and Viceira(1999)]{Campbell_Viceira99}
J.Y. Campbell and L.M. Viceira.
\newblock Consumption and portfolio decisions when expected returns are time
  varying.
\newblock {\em Quarterly Journal of Economics}, 114:\penalty0 433--495, 1999.

\bibitem[Campbell and Viceira(2001)]{Campbell_Viceira99b}
J.Y. Campbell and L.M. Viceira.
\newblock Who should buy long-term bonds?
\newblock {\em American Economic Review}, 91:\penalty0 99--127, 2001.

\bibitem[Chabi-Yo et~al.(2005)Chabi-Yo, Garcia, and Renault]{Chabi-yo_etal05}
Fousseni Chabi-Yo, René Garcia, and Eric Renault.
\newblock State dependence in fundamentals and preferences explains
  risk-aversion puzzle.
\newblock Working Papers 05-9, Bank of Canada, 2005.
\newblock available at http://ideas.repec.org/p/bca/bocawp/05-9.html.

\bibitem[Constantinides(1990)]{Constantinides90}
G.~M. Constantinides.
\newblock Habit formation: A resolution of the equity premium puzzle.
\newblock {\em Journal of Political Economy}, 98:\penalty0 519--543, 1990.

\bibitem[Cox and Huang(1989)]{Cox_etal89}
John~C. Cox and Chi-fu Huang.
\newblock Optimal consumption and portfolio policies when asset prices follow a
  diffusion process.
\newblock {\em J. Econom. Theory}, 49\penalty0 (1):\penalty0 33--83, 1989.
\newblock ISSN 0022-0531.

\bibitem[Cvitani{\'c} et~al.(2003)Cvitani{\'c}, Goukasian, and
  Zapatero]{Cvitanic_etal03}
Jak{\v{s}}a Cvitani{\'c}, Levon Goukasian, and Fernando Zapatero.
\newblock Monte {C}arlo computation of optimal portfolios in complete markets.
\newblock {\em J. Econom. Dynam. Control}, 27\penalty0 (6):\penalty0 971--986,
  2003.
\newblock ISSN 0165-1889.
\newblock High-performance computing for financial planning (Ischia, 1999).

\bibitem[Dammon et~al.(2001)Dammon, Spatt, and Zhang]{Dammon_etal01}
Robert~M. Dammon, Chester~S. Spatt, and Harold~H. Zhang.
\newblock Optimal consumption and investment with capital gains taxes.
\newblock {\em Review of Financial Studies}, 14:\penalty0 583--616, 2001.

\bibitem[Dangl and Wirl(2004)]{Dangl_Wirl04}
Thomas Dangl and Franz Wirl.
\newblock Investment under uncertainty: calculating the value function when the
  {B}ellman equation cannot be solved analytically.
\newblock {\em J. Econom. Dynam. Control}, 28\penalty0 (7):\penalty0
  1437--1460, 2004.
\newblock ISSN 0165-1889.
\newblock Mathematical programming.

\bibitem[Dantine et~al.(2004)Dantine, Donaldson, Giannikos, and
  Guirguis]{Danthine_etal04}
Jean-Pierre Dantine, John~B. Donaldson, Christos Giannikos, and Hany Guirguis.
\newblock On the consequences of state dependent preferences for the pricing of
  financial assets.
\newblock {\em Finance Research Letters}, 1:\penalty0 143--153, 2004.

\bibitem[Davis and Norman(1990)]{Davis_Norman90}
M.~H.~A. Davis and A.~R. Norman.
\newblock Portfolio selection with transaction costs.
\newblock {\em Math. Oper. Res.}, 15\penalty0 (4):\penalty0 676--713, 1990.
\newblock ISSN 0364-765X.

\bibitem[Detemple et~al.(2003)Detemple, Garcia, and
  Rindisbacher]{Detemple_etal03}
J{\'e}r{\^o}me~B. Detemple, Ren{\'e} Garcia, and Marcel Rindisbacher.
\newblock A monte carlo method for optimal portfolios.
\newblock {\em Journal of Finance}, 58\penalty0 (1):\penalty0 401--446, 2003.

\bibitem[Duffie and Epstein(1992)]{Duffie_and_Epstein92a}
D.~Duffie and L.~(Appendix C with C.~Skiadas) Epstein.
\newblock Stochastic differential utility.
\newblock {\em Econometrica}, 60:\penalty0 353--394, 1992.

\bibitem[Epstein and Zin(1989)]{Epstein_and_Zin89}
L.~Epstein and S.~Zin.
\newblock Substitution, risk advertion and the temporal behavior of asset
  returns: A theoretical framework.
\newblock {\em Econometrica}, 57:\penalty0 937--969, 1989.

\bibitem[Hindy and Huang(1993)]{Hindy_and_Huang93}
A.~Hindy and C.~F. Huang.
\newblock Optimal consumption and portfolio rules with durability and local
  substitution.
\newblock {\em Econometrica}, 61:\penalty0 85--122, 1993.

\bibitem[Hindy et~al.(1997)Hindy, Huang, and Zhu]{Hindy_etal97}
Ayman Hindy, Chi-fu Huang, and Steven~H. Zhu.
\newblock Numerical analysis of a free-boundary singular control problem in
  financial economics.
\newblock {\em J. Econom. Dynam. Control}, 21\penalty0 (2-3):\penalty0
  297--327, 1997.
\newblock ISSN 0165-1889.

\bibitem[Hindya et~al.(1997)Hindya, fu~Huanga, and Zhu]{Hindy_etal97b}
Ayman Hindya, Chi fu~Huanga, and Steven~H. Zhu.
\newblock Optimal consumption and portfolio rules with durability and habit
  formation.
\newblock {\em Journal of Economic Dynamics and Control}, 21\penalty0
  (2--3):\penalty0 525--550, 1997.

\bibitem[Jackwerth(2000)]{Jackwerth00}
J.~C. Jackwerth.
\newblock Recovering risk aversion from option prices and realized returns.
\newblock {\em Review of Financial Studies}, 13\penalty0 (2):\penalty0
  433--451, 2000.

\bibitem[Karatzas et~al.(1987)Karatzas, Lehoczky, and Shreve]{Karatzas_etal87}
Ioannis Karatzas, John~P. Lehoczky, and Steven~E. Shreve.
\newblock Optimal portfolio and consumption decisions for a ``small investor''
  on a finite horizon.
\newblock {\em SIAM J. Control Optim.}, 25\penalty0 (6):\penalty0 1557--1586,
  1987.
\newblock ISSN 0363-0129.

\bibitem[Karatzas and Shreve(1998)]{Karatzas98}
Ioannis Karatzas and Steven~E. Shreve.
\newblock {\em Methods of Mathematical Finance}, volume~39 of {\em Applications
  of Mathematics}.
\newblock Springer-Verlag, New York, 1998.

\bibitem[Karni(1993{\natexlab{a}})]{Karni93b}
Edi Karni.
\newblock A definition of subjective probabilities with state-dependent
  preferences.
\newblock {\em Econometrica}, 61\penalty0 (1):\penalty0 187--198,
  1993{\natexlab{a}}.
\newblock ISSN 0012-9682.

\bibitem[Karni(1993{\natexlab{b}})]{Karni93}
Edi Karni.
\newblock Subjective expected utility theory with state-dependent preferences.
\newblock {\em J. Econom. Theory}, 60\penalty0 (2):\penalty0 428--438,
  1993{\natexlab{b}}.
\newblock ISSN 0022-0531.

\bibitem[Kim and Omberg(1996)]{Kim_Omberg96}
Tong~Suk Kim and Edward Omberg.
\newblock Dynamic nonmyopic portfolio behavior.
\newblock {\em Review of Financial Studies}, 9:\penalty0 141--161, 1996.

\bibitem[Kunita(1990{\natexlab{a}})]{Kunita90}
Hiroshi Kunita.
\newblock {\em Stochastic flows and stochastic differential equations},
  volume~24 of {\em Cambridge Studies in Advanced Mathematics}.
\newblock Cambridge University Press, Cambridge, 1990{\natexlab{a}}.
\newblock ISBN 0-521-35050-6.

\bibitem[Kunita(1990{\natexlab{b}})]{hK90}
Hiroshi Kunita.
\newblock {\em Stochastic flows and stochastic differential equations},
  volume~24 of {\em Cambridge studies in advanced mathematics}.
\newblock Cambridge University Press, Cambridge, New York, Melbourne,
  1990{\natexlab{b}}.

\bibitem[Lazrak and Quenez(2003)]{Lazrak03}
Ali Lazrak and Marie~Claire Quenez.
\newblock A generalized stochastic differential utility.
\newblock {\em Math. Oper. Res.}, 28\penalty0 (1):\penalty0 154--180, 2003.
\newblock ISSN 0364-765X.

\bibitem[Lioui and Poncet(2001)]{Abramham_Poncet01}
Abraham Lioui and Patrice Poncet.
\newblock On optimal portfolio choice under stochastic interest rates.
\newblock {\em J. Econom. Dynam. Control}, 25\penalty0 (11):\penalty0
  1841--1865, 2001.
\newblock ISSN 0165-1889.

\bibitem[Londo{\~n}o(2004)]{Londono04}
Jaime~A. Londo{\~n}o.
\newblock State tameness: A new approach for credit constrains.
\newblock {\em Electronic Communications in Probability}, 9:\penalty0 1--13,
  January 2004.

\bibitem[Londoño(2006{\natexlab{a}})]{Londono2006}
Jaime~A. Londoño.
\newblock Dynamic state dependent utilities.
\newblock Working paper, 2006{\natexlab{a}}.

\bibitem[Londoño(2006{\natexlab{b}})]{Londono2005a}
Jaime~A. Londoño.
\newblock Dynamic state tameness.
\newblock Available at
  http://math.ucr.edu/$\sim$jlondono/papers/dynamic\_tameness.pdf,
  2006{\natexlab{b}}.

\bibitem[Magill and Constantinides(1976)]{Magill_Constantinides76}
Michael J.~P. Magill and George~M. Constantinides.
\newblock Portfolio selection with transactions costs.
\newblock {\em J. Econom. Theory}, 13\penalty0 (2):\penalty0 245--263, 1976.
\newblock ISSN 0022-0531.

\bibitem[Mehra and Prescott(1985)]{Mehra_and_Prescott85}
R.~Mehra and E.~Prescott.
\newblock The equity premium: A puzzle.
\newblock {\em Journal of Monetary Economics}, 15:\penalty0 145--161, 1985.

\bibitem[Melino and Yang(2003)]{Melino03}
Angelo Melino and Alan~X. Yang.
\newblock State dependent preferences can explain the equity premium puzzle.
\newblock {\em Review of Economic Dynamics}, 6:\penalty0 806--830, 2003.

\bibitem[Merton(1969)]{Merton69}
Robert~C. Merton.
\newblock Lifetime portfolio selection under uncertainty: the continuos time
  case.
\newblock {\em Review of Economics and Statistics}, 51\penalty0 (3):\penalty0
  247--257, 1969.

\bibitem[Merton(1971)]{Merton71}
Robert~C. Merton.
\newblock Optimum consumption and portfolio rules in a continuous-time model.
\newblock {\em J. Econom. Theory}, 3\penalty0 (4):\penalty0 373--413, 1971.
\newblock ISSN 0022-0531.

\bibitem[Ocone and Karatzas(1991)]{Ocone_Karatzas91}
Daniel~L. Ocone and Ioannis Karatzas.
\newblock A generalized {C}lark representation formula, with application to
  optimal portfolios.
\newblock {\em Stochastics Stochastics Rep.}, 34\penalty0 (3-4):\penalty0
  187--220, 1991.
\newblock ISSN 1045-1129.

\bibitem[Shreve and Soner(1994)]{Shreve_Soner94}
S.~E. Shreve and H.~M. Soner.
\newblock Optimal investment and consumption with transaction costs.
\newblock {\em Ann. Appl. Probab.}, 4\penalty0 (3):\penalty0 609--692, 1994.
\newblock ISSN 1050-5164.

\bibitem[Watchter(2002)]{Wachter02}
Jessica Watchter.
\newblock Portfolio and consumption decisions under mean-reverting returns: An
  exact solution for complete markets.
\newblock {\em Journal of Financial and Quantitative Analysis}, 37\penalty0
  (1):\penalty0 63--91, 2002.

\bibitem[Weil(1989)]{Weil89}
O.~Weil.
\newblock The equity premium puzzle and the risk-free rate puzzle.
\newblock {\em Journal of Monetary Economics}, 24\penalty0 (3):\penalty0
  401--421, 1989.

\end{thebibliography}
\end{document}